
\input amstex
\documentstyle{amsppt}
\magnification=1100
\pagewidth{135mm}
\pageheight{190mm}

\topmatter
\title
A remark on positively curved manifolds of
dimensions 7 and 13
\endtitle
\author
Iskander A. TAIMANOV
\endauthor
\address
Institute of Mathematics,
630090 Novosibirsk, Russia
\endaddress
\email
taimanov\@math.nsk.su
\endemail
\endtopmatter
\leftheadtext{}
\rightheadtext{}
\document
\NoBlackBoxes

In the present paper we construct totally geodesic embeddings of some
7-dimensional manifolds into 13-dimensional manifolds with positive
sectional curvature
and explain the strange coincidence of pinching constants of the normally
homogeneous Berger space (\cite{Be}) and the homogeneous Aloff-Wallach
space $N_{1,1,-1/2}$ (\cite{AW}).
This constant is equal to $\frac{16}{29 \cdot 37}$
as it was established by Heintze for the Berger space (\cite{H}) and by
Huang (\cite{Hu}) for the Aloff-Wallach space.
Moreover these totally geodesic embeddings shed light on a relation
of the well-known 7-dimensional manifolds constructed by Aloff,
Wallach and Eschenburg (\cite{E1}) to the series of 13-dimensional
positively curved manifolds founded recently by Bazaikin (\cite{Ba}).

We identify the Lie groups $U(n)$
with the groups formed by $(n \times n)$-matrices $A$ such that
$$
A \cdot I_n \cdot A^* = I_n
$$
where $I_n$ is the unit $(n \times n)$-matrix and $A^*_{ij} = {\bar A}_{ji}$.
We consider the Lie groups $SU(n)$ as the subgroups of $U(n)$
formed by matrices with $\det =1$.
We mean by $Sp(2)$ a subgroup of $SU(4)$ formed by matrices $A$
which satisfy the following equality
$$
A \cdot
\left(
\matrix 0 & I_2 \\ -I_2 & 0 \endmatrix
\right)
\cdot A^t =
\left(
\matrix 0 & I_2 \\ -I_2 & 0 \endmatrix
\right)
$$
where $A^t_{ij} = A_{ji}$.
Moreover we will consider $Sp(2)$ as a subgroup of $SU(5)$ formed by
matrices of the form
$$
\left(
\matrix
A & 0 \\
0 & 1
\endmatrix
\right) \in SU(5), \ \ A \in Sp(2) \subset SU(4).
$$

We consider the Lie algebras of these groups as realized by matrices
in usual manner.

We denote by $T^1,T^1_p$ and $T_{(0)}$
the subgroups of $U(5)$ formed by diagonal matrices
\linebreak
${\textstyle diag}(z,z,z,z,z^{-4})$, ${\textstyle diag}(z^p,z^p,z^p,z^p,z)$,
and ${\textstyle diag}(z,z,z,z,1)$
respectively where $|z|=1$ and $p$ is a positive integer number.

First we remind the known examples (we consider only simply connected
manifolds).

1) 7-dimensional manifolds.

1.1) 7-dimensional Berger space (\cite{Be}).

This space is isometric to
a factor-space $Sp(2)/SU(2)$ where $Sp(2)$ is endowed with a standard
biinvariant metric and an embedding $SU(2) \subset Sp(2)$ is a nonstandard
one. We will not discuss this space and
consider it as an exceptional one.

1.2) Aloff-Wallach spaces (\cite{AW}).

Let $T_{k,l}$ be the subgroup of $SU(3)$ formed by diagonal matrices
$(z^k,z^l,z^{-(k+l)})$. We consider the subgroup $G_1=U(2)$ of $SU(3)$
given by
$$
G_1 = \left( \matrix A & 0 \\ 0 & \det A^{-1} \endmatrix
\right), \ \ \
A \in U(2)
$$
and denote by $g_1$ the Lie algebra of $G_1$. We denote by $f_{k,l}$
the Lie algebra of $T_{k,l}$ . One can see that $f_{k,l}$ is generated
by a diagonal matrix ${\textstyle diag}(2\pi \sqrt{-1} k, 2\pi\sqrt{-1}
l, -2\pi\sqrt{-1}(k+l))$.  Let denote by $\langle,\rangle_0$ the
Killing biinvariant metric on $SU(3)$.  Then one can consider a
homogeneous metric, on factor-space $N_{k,l,t}=SU(3)/T_{k,l}$,
generated by the metric
$$
\langle x,y \rangle = (1+t) \langle x_1,y_1
\rangle_0 + \langle x_2,y_2 \rangle_0
\eqno{(1)}
$$
where $x_i,y_i \in
V_i$ and $f_{k,l}^{\perp} = V_1 \oplus V_2$ is an orthogonal
decomposition, $V_1 = f_{k,l}^{\perp} \cap g_1$ and $V_2 =
g_1^{\perp}$.

Aloff and Wallach
had shown that if $k$ and $l$ have the same sign
and $-1 < t <0$ then these manifolds $N_{k,l,t}$ are positively curved.

1.3) Eschenburg spaces (\cite{E1,E2}).

These spaces are generalizations of the previous ones and have form
$T_{k,l} \backslash U(3) / T_{p,q}$. For suitable  integers $k,l,p$ and $q$
and a metric on $U(3)$ these manifolds would have positive curvature.
We will not dwell on them and only
notice that most of them are not homeomorphic
to homogeneous manifolds and these manifolds were first examples of such
kind.

2) 13-dimensional manifolds.

2.1) 13-dimensional Berger space (\cite{Be}).

This manifold $B^{13}$ is a factor-space $SU(5) / (Sp(2) \times T^1)$ where
$SU(5)$ is endowed with the Killing biinvariant metric.

2.2) Bazaikin spaces (\cite{Ba}).

These spaces have form
$ S_{\bar p} \backslash U(5) / (Sp(2) \times T_{(0)})$
where $S_{\bar p}$ is formed by diagonal matrices
${\textstyle diag}(z^{p_1},z^{p_2},z^{p_3},z^{p_4},z^{p_5})$
where $|z|=1$.
For suitable tuples of integers $\bar p$ and a metric on $U(5)$ these
manifolds would have positive curvature.
We also mention that for these examples a metric on $U(5)$ is taken to be
left-invariant under $U(5)$-action and right-invariant under
$U(4) \times U(1)$-action where subgroup $U(4) \times U(1)$ formed by
diagonal-block matrices with $4 \times 4$- and $1 \times 1$- blocks.

Now we proceed with our main construction.

Put
$$
S =
\left(
\matrix
0 & 0 & 1 & 0 & 0 \\
0 & 0 & 0 & 1 & 0 \\
-1& 0 & 0 & 0 & 0 \\
0 &-1 & 0 & 0 & 0 \\
0 & 0 & 0 & 0 & \sqrt{-1}
\endmatrix
\right)
$$
and define the following mapping
$$ \sigma : U(5) \rightarrow U(5) \ \ \ : \ \ \
A \rightarrow S \cdot A \cdot S^{-1}.
$$

The following two propositions are evident.

\proclaim{Proposition 1}
$\sigma(G) = G$ for $G = SU(5), Sp(2), T^1, T^1_p, T_{(0)}$.
\endproclaim

\proclaim{Proposition 2}
$\sigma^2 $ is an identical mapping, i.e., $\sigma$ is an involution.
\endproclaim

Since metric on $SU(5)$ is biinvariant for the Berger space and
the metric on
is $U(5)$-left-invariant and $U(4) \times U(1)$-right-invariant for
the Bazaikin spaces $T^1_p \backslash U(5) / (Sp(2) \times T_{(0)})$,
the following Proposition holds.

\proclaim{Proposition 3}
The involution $\sigma$ induces isometric involutions on the spaces
$B^{13}$ and
$T^1_p \backslash U(5) / (Sp(2) \times T_{(0)}) (p >0)$.
\endproclaim

First we consider the action of this involution on the
Berger space $B^{13}$.

\proclaim{Theorem 1}
Let $W^7$ be a submanifold of $B^{13}$
which contains the point $E = 1 \cdot (Sp(2) \times T^1)
\in B^{13}$
where $1$ is the unit of $SU(5)$ and which is
formed by fixed points of involution $\sigma : B^{13} \rightarrow B^{13}$.
Then the manifold $W^7$ is a totally geodesic submanifold which is
isometric to the Aloff-Wallach space
$N_{1,1,-1/2}$ and minimal and maximal values of the sectional curvature
of $B^{13}$ are attained on two-planes tangent to $W^7$.
\endproclaim

Proof of Theorem 1.

First notice that, since $W^7$ is a component of  the set formed by
fixed points of involution, this embedding $W^7 \rightarrow B^{13}$ is
a totally geodesic one.

Now let compute the dimension of $W^7$ and find generators of the tangent
space of $W^7$ at the point $E$.

We use notations from the paper of Heintze (\cite{H}) who denoted by
$H_i ( 1 \leq i \leq 11)$ a set of orthonormal  vectors which form a
basis for tangent space to $Sp(2) \times T^1$ and denoted by $M_j
(1\leq j \leq 13)$ basic orthonormal vectors of its orthogonal
complement.

Let $E_{kl}$ be a $(5\times 5)$-matrix
$(\delta_{ak}\delta_{bl})_{(1 \leq a,b \leq 5)}$.
Then let introduce $Q_{kl} = E_{kl}-E_{lk}, R_{kl}= \sqrt{-1}(E_{kl}+E_{lk})$
and $P_k = \sqrt{-1}(E_{kk} - E_{55})$.

Heintze used the following basis :
$$
M_j = \sqrt{2} Q_{j5},\ \  M_{j+4} = \sqrt{2} R_{j5} , \ \ j=1,2,3,4,
$$
$$
M_9 = Q_{12}-Q_{34} ,\ \  M_{10} = Q_{14}-Q_{23},
$$
$$
M_{11} = R_{12}+R_{34},\ \  M_{12} = R_{14} - R_{23} ,
\ \ M_{13} = P_1-P_2+P_3-P_4.
$$

One can derive by direct computations that the space $V$ generated by these
vectors splits into two pairwise orthogonal subspaces $V^+$ and $V^-$ such
that $\sigma|_{V^{\pm}} = \pm  1$.

The orthonormal bases of these subspaces are :

1) for $V^+$ :
$$
\frac{M_1+M_7}{\sqrt{2}},\frac{M_2-M_8}{\sqrt{2}},
\frac{M_3-M_5}{\sqrt{2}},\frac{M_4-M_6}{\sqrt{2}},
M_{11},M_{12},M_{13} ;
$$

2) for $V^-$ :
$$
\frac{M_1-M_7}{\sqrt{2}},\frac{M_2-M_8}{\sqrt{2}},
\frac{M_3+M_5}{\sqrt{2}},\frac{M_4+M_6}{\sqrt{2}},
M_9,M_{10}.
$$

Since $\dim V^+ =7$ and the submanifold $W^7$ is homogeneous,
$$
 \dim W^7 =7.
$$

Let show that this submanifold is isometric to the space $N_{1,1,-1/2}$.

We introduce another action given by
$$
\rho: SU(5) \rightarrow  SU(5)
\ \ \ \  :
\ \ \
A \rightarrow R \cdot A \cdot R^{-1}
$$
where
$$
R =
\frac{1}{\sqrt{2}}
\left(
\matrix
1 & 0 & \sqrt{-1} & 0 & 0 \\
0 & 1 & 0 & \sqrt{-1} & 0 \\
\sqrt{-1} & 0 & 1 & 0 & 0 \\
0 & \sqrt{-1} & 0 & 1 & 0 \\
0 & 0 & 0 & 0 & \sqrt{2}
\endmatrix
\right).
$$

Since $R \in Sp(2) \subset  SU(5)$,
this action is an isometry.

The action $\rho$ generates an action on the Lie algebra $su(5)$
which we denote also by $\rho$.
We compute the action of $\rho$ only on $V^+$ :
$$
\rho(\frac{M_1+M_7}{\sqrt{2}}) = M_7,\ \ \
\rho(\frac{M_2+M_8}{\sqrt{2}}) = M_8,
$$
$$
\rho(\frac{M_3-M_5}{\sqrt{2}}) = M_3,\ \
\rho(\frac{M_4-M_6}{\sqrt{2}}) = M_4,
\eqno{(2)}
$$
$$
\rho(M_{11}) = M_{11},\ \  \rho(M_{12}) = M_9 ,\ \
\rho(M_{13}) = M_{13}.
$$

Moreover we have
$$
M_9 = -2 Q_{34} + \rho(H_{10}),\
M_{11} = 2 R_{34} + \rho(H_8),\
M_{13} = 2(P_3-P_4) + \rho(H_5)
\eqno{(3)}
$$
where $H_j$ are unit basic orthonormal
vectors from the tangent space to $Sp(2) \times T^1$
(see \cite{H}).

One can see that
$(M_3,M_4,M_7,M_8, \sqrt{2} Q_{34}, \sqrt{2} R_{34}, \sqrt{2}(P_3-P_4))$
form an orthonormal (up to multiplication by constant) basis in
a subalgebra $V_1$ (see (1))
for $k=l=1$ and the group $SU(3)$ given by
$$
\tilde G =
\left(
\matrix
I_2 & 0 \\
0 & A
\endmatrix
\right),
\ \ \ \
A \in SU(3).
\eqno{(4)}
$$

But , since (2) and (3), we have that
$(M_3,M_4,M_7,M_8, 2Q_{34}, 2R_{34},2(P_3-P_4))$
form an orthonormal basis in the tangent space of $\rho(W)$ at
the point $E$.
This coincides with (1) for $t=-\frac{1}{2}$.

Since $W^7$ is a homogeneous manifold and
the homogeneous manifold $N_{1,1,-1/2}$
is simply connected, there exists a finite isometric covering
$$
N_{1,1,-1/2} \rightarrow W^7.
$$

Let prove that this covering is a diffeomorphism.

{}From formulas (2) and (3) one can see that $W^7$ is formed by
$Sp(2) \times T^1$-orbits
of elements $g$ such that $\rho(g) \in \tilde G$ (see (4)).
But by direct computations one can derive that

1) if $\rho(h) \in \tilde G$ and $h \in Sp(2)$ then $h = 1 \in Sp(2)$ ;

2) if
$$
h = \left(
\matrix a & b  \\ c & d \endmatrix
\right) \in Sp(2)
$$
and $\rho(h)$ is a diagonal matrix then $a=d, b=-c , a^2+b^2 = I_2$ and
$\rho(h)$ has the form
${\textstyle diag}(\lambda,\mu,\lambda^{-1},\mu^{-1},1)$.

It follows now that orbits of two elements $g_1,g_2 \in\tilde G$ coincide
if and only if $g_1g_2^{-1} \in T_{1,1} \subset SU(3)$ and we conclude that
$$
W^7 = N_{1,1,-1/2}.
$$

We are left to prove that the curvature of $B^{13}$ attains its minimal and
maximal values on two-planes which are tangent to $W^7$.

By using of the formula from Lemma 2 of \cite{H} one compute the sectional
curvature of 2-plane generated by vectors $X$ and $Y$ :
$$
K(X,Y) = \frac{29}{4}
$$
for $X = \frac{M_1+M_7}{\sqrt{2}}, Y=\frac{M_3-M_5}{\sqrt{2}}$ and
it was proved in \cite{H} that this value is the maximum of curvature of
the space $B^{13}$.

Let take a matrix
$Q={\textstyle diag}(1,-\sqrt{-1},1,\sqrt{-1},1) \in Sp(2) \subset SU(5)$.
Since $Q \in Sp(2)$, the action $\xi: X \rightarrow Q \cdot X \cdot Q^{-1}$
generates an isometry of $B^{13}$ (see \cite{H}).
Let take
$$
X = \sqrt{\frac{12}{37}}(\frac{M_1+M_7}{\sqrt{2}} + \frac{M_2+M_8}{\sqrt{2}})
+ \sqrt{\frac{13}{37}}M_{11},
$$
$$
Y = - \sqrt{\frac{12}{37}}(\frac{M_3-M_5}{\sqrt{2}} -
\frac{M_4-M_6}{\sqrt{2}}) - \sqrt{\frac{13}{37}} M_{12}.
$$
One can compute
$$
\xi^{-1}(X) = \sqrt{\frac{12}{37}}(\frac{M_1+M_7}{\sqrt{2}} +
\frac{M_4+M_6}{\sqrt{2}}) + \sqrt{\frac{13}{37}}M_9,
$$
$$
\xi^{-1}(Y) = -\sqrt{\frac{12}{37}}(\frac{M_3-M_5}{\sqrt{2}} -
\frac{M_2-M_8}{\sqrt{2}}) + \sqrt{\frac{13}{37}}M_{10},
$$
apply to these vectors  Lemma 2 of \cite{H} and
derive that
$$
K(X,Y)=K(\xi^{-1}(X),\xi^{-1}(Y))=\frac{4}{37}.
$$
It was proved in \cite{H} that this value is the minimum of the
curvature of $B^{13}$.

One can conclude that the pinching constants $K_{\min}/K_{\max}$
for $B^{13}$ and $W^7=N_{1,1,-1/2}$ coincide and are equal to
$16 / 29 \cdot 37$.

Theorem 1 is established.

Of course, this explanation can be simplified by replacing the involution
$\sigma$ by another one:
$$
A \rightarrow \Sigma \cdot A \cdot \Sigma^{-1}, \ \
\Sigma = R \cdot S \cdot R^{-1}.
$$
Nevertheless we prefer to describe this proof subsequently as it was
originally obtained.

By using analogous conversations one can prove the following
theorem.

\proclaim{Theorem 2}
Let $M^{13}_p$ be the Bazaikin space of the form
$T^1_p \backslash U(5) / Sp(2) \times T_{(0)}$.
Then the fixed-point set, of involution
$\sigma: M^{13}_p \rightarrow M^{13}_p$, which contains the point
$E = T^1_p \cdot 1 \cdot Sp(2) \times T_{(0)}$,
is diffeomorphic to the space
$$
W^7_p =
\left(
\matrix
z^{p} & 0 & 0 \\
0 & z^{p} & 0 \\
0 & 0 & z
\endmatrix
\right)
{\big \backslash U(3) /}
\left(
\matrix
{\bar w} & 0 & 0 \\
0 & {\bar w} & 0 \\
0 & 0 & 1
\endmatrix
\right), \ \ \ \ |z|=|w|=1.
$$
\endproclaim

One can easy see that $W^7_p$ is a totally geodesic submanifold of $M^{13}_p$.
The spaces $M^{13}_p$ are nonhomogeneous for $p \geq 2$ and thus the problem of
comparing pinchings of $W^7_p$ and $M^{13}_p$ is not reduced to local
computations as it was done in the proof of Theorem 1.

{\bf Remark.} These spaces $W^7_p$ are not presented directly in this form in
\cite{E1,E2} and probably some of these examples were not known before.
We notice that they are also of the biquotient form which was
introduced by Gromoll and Meyer (\cite{GM}). As they are totally geodesic
submanifolds of positively curved space they have positive
sectional curvature.

The following question looks interesting.

\proclaim{Question}
Does there exist a correspondence of 7-dimensional Aloff-Wallach and
Eschenburg spaces to 13-dimensional Berger and Bazaikin spaces
which is realized by totally geodesic embeddings ?

If such correspondence exists only for some subfamilies what are these
subfamilies ?

If such correspondence exists is it realized by pinching-essential
embeddings (i.e., embeddings with the same pinching constants of
manifolds and submanifolds) as in the case of Theorem 1 ?
\endproclaim

Let consider the topological properties of manifolds $W^7$ and $B^{13}$.

Put
$$
\hat G (\approx SU(2)) = \left(\matrix
A & 0 \\ 0 & 1 \endmatrix\right) \subset SU(3), \ \
A \in SU(2).
$$
One can see that $SU(3) / T_{1,1} = W^7$ and $SU(3)/{\hat G} = S^5$.
The group $T_{1,1}$ acts on
$S^5 = \{z_1^2+z_2^2+z_3^2=1 | z_i \in {\bold C}\}$ by multiplications :
$$
(z_1,z_2,z_3) \rightarrow (\lambda^{-2} z_1, \lambda^{-2} z_2,
\lambda^{-2} z_3)
$$
where ${\textstyle diag}(\lambda,\lambda,\lambda^{-2}) \in T_{1,1}$.
Moreover the actions of $T_{1,1}$ and $\hat G$, on $SU(3)$, commute.

Let us consider the fiber bundle
$$
SU(3) \rightarrow {\bold C}P^2.
\eqno{(5)}
$$
Its fiber is diffeomorphic to $U(2)$ which one can represent in the
following manner. Put
$$
\hat Q = SU(2) \times  ({\bold R}/2\pi{\bold Z}).
$$
We denote by $\hat Q_1$ the factor-space
of $\hat Q$ under the following
${\bold Z}_2$-action:
$$
(X,t) \rightarrow (-X, t+\pi), \ \
X \in SU(2).
$$
This factor-space is diffeomorphic to the fiber of bundle (5).
It is fibered over $S^1$ in the usual manner
$$
(X,t) \rightarrow t \in S^1 = {\bold R}/\pi{\bold Z}.
$$
In these terms the action of $T_{1,1}$ on fiber bundle (5)
has the form
$$
(X,t) \rightarrow (\exp{(\sqrt{-1}\pi \phi)} \cdot X, t+\phi),
$$
$$
{\textstyle diag}(\exp{(\sqrt{-1}\pi \phi)},\exp{(\sqrt{-1}\pi \phi)},
\exp{(-2\sqrt{-1}\pi \phi)}) \in T_{1,1}.
$$

Now one can derive that

1) $S^5 / T_{1,1} = SU(3) / {\hat G} \times T_{1,1} = {\bold C}P^2$ ;

2) these actions generate a fiber bundle
$$
W^7 = SU(3) / T_{1,1} \overset{{\bold R}P^3}\to{\longrightarrow}
{\bold C}P^2 ;
\eqno{(6)}
$$

3) it follows from computations of cohomology
groups of $W^7$ (see \cite{E1})
that the transgression $d_4$ in
the spectral sequence of  fiber bundle
(6) is given by
$$
d_4 : E^{0,3}_4  = {\bold Z} \overset{\times
3}\to{\longrightarrow} E^{4,0}_4  = {\bold Z},
\eqno{(7)}
$$
and
$$
H^4(W^7)={\bold Z}_3.
\eqno{(8)}
$$

Put
$$
\bar G (\approx SU(4)) =
\left(\matrix A & 0 \\ 0 & 1
\endmatrix\right) \subset SU(5), \ \
A \in SU(4).
$$
One can see that $SU(5) / Sp(2) \times T^1 = B^{13}$ and
$SU(5)/{\bar G} = S^9$.
The group $T^1$ acts on
$S^9=\{z_1^2 + \dots + z_5^2 = 1 | z_i \in {\bold C}\}$
by multiplications :
$$
(z_1,\dots,z_5) \rightarrow (\lambda^{-4} z_1,\dots,\lambda^{-4} z_5)
$$
where ${\textstyle diag}(\lambda,\lambda,\lambda,\lambda,\lambda^{-4})
\in T^1$. Moreover the actions of $T^1$ and ${\bar G}$, on $SU(5)$,
commute.

Let us consider the fiber bundle
$$
SU(5)/Sp(2) \rightarrow {\bold C}P^4.
\eqno{(9)}
$$
Put
$$
\bar Q = SU(4)/Sp(2) \times  ({\bold R}/\pi{\bold Z}).
$$
We denote by $\bar Q_1$ the factor-space
of $\bar Q$ under the following
${\bold Z}_2$-action:
$$
(X,t) \rightarrow (\sqrt{-1}X, t+\frac{\pi}{2}), \ \
X \in SU(4)/Sp(2).
$$
This factor-space is diffeomorphic to the fiber of bundle (9).
It is fibered over $S^1$ in the usual manner
$$
(X,t) \rightarrow t \in S^1 = {\bold R}/\frac{\pi}{2}{\bold Z}.
$$
In these terms the action of $T^1$ on fiber bundle (9)
has the form
$$
(X,t) \rightarrow (\exp{(\sqrt{-1}\pi \phi)} \cdot X, t+\phi),
$$
$$
{\textstyle diag}(\exp{(\sqrt{-1}\pi \phi)}, \dots,
\exp{(\sqrt{-1}\pi \phi)}, \exp{(-4\sqrt{-1}\pi \phi)}) \in T^1.
$$

Now one derive that

1) $S^5 / T^1 = SU(5) / {\bar G} \times T^1 = {\bold C}P^4$ ;

2) these actions generate a fiber bundle
$$
B^{13} \overset{{\bold R}P^5}\to{\longrightarrow}
{\bold C}P^4 ; \eqno{(10)}
$$

3) it follows from computations of cohomology groups of $B^{13}$ (see
\cite{Ba}) that the transgression $d_6$
in the spectral sequence of fiber bundle (10) is given by
$$
d_6: E^{0,5}_6 = {\bold Z} \overset{\times
5}\to{\longrightarrow} E^{6,0}_6 =  {\bold Z}.
\eqno{(11)}
$$
and
$$
H^6(B^{13})={\bold Z}_5.
\eqno{(12)}
$$

The similarity of formulas (5-8) for $W^7$ and formulas (9-12)
for $B^{13}$ and Theorem 1 give us a reason to pose the following
question.

\proclaim{Question}
Is it true that for every positive integer $k$ there exist a space
$\Gamma_k$ such that

1) there exists a fiber bundle
$$
\Gamma_k \overset{{\bold R}P^{2k+1}}\to{\longrightarrow}
{\bold C}P^{2k} ; \eqno{(13)} $$

2) the transgression $d_{2k+2}$
in the spectral sequence of (13) is given by
$$
d_{2k+2} : H^{2k+1}(S^{2k+1}) \overset{\times (2k+1)}
\to{\longrightarrow} H^{2k+2}({\bold C}P^{2k}),
\eqno{(14)}
$$
and
$$
H^{2k+2}(\Gamma_k)={\bold Z}_{2k+1} \ \ ;
\eqno{(15)}
$$

3) a manifold $\Gamma_k$ has positive sectional curvature ;

4) $\Gamma_1 = W^7$ and $\Gamma_2=B^{13}$ ?
\endproclaim

One can pose more rigorous conjecture by adding

{\sl  5) the spaces $\Gamma_k$ form a tower
$$
\Gamma_1 \rightarrow \Gamma_2
\rightarrow \dots \rightarrow
\Gamma_n \rightarrow \Gamma_{n+1} \rightarrow \dots
$$
of pinching-essential totally geodesic embeddings and, thus,
the pinching constants of $\Gamma_k$ are equal to $\frac{16}{29 \cdot
37}$.}

If such tower exists one
can expect that its properties are similar to the properties
of ${\bold C}P^n$- or ${\bold H}P^n$-towers.

{\bf Final remarks.}

1)  Let consider topological properties of the spaces
$W^7_p$ and $M^{13}_p$.

Since a left-side multiplication  transforms orbits
under right-side action into orbits under this action,
the group $\tilde T_p = {\textstyle diag}(z^p,z^p,z)$ acts on the space
$U(3)/SU(2) \cdot {\textstyle diag}(\bar w,\bar w,1) $
\newline
$= S^5$ (where $|z|=|w|=1$) and the group $T^1_p$ acts on the space
$U(5)/SU(4) \cdot T_{(0)} = S^9$.

These actions are not free.

One can see that the elements ${\textstyle diag}(z^p,z^p,z)$ for
$z^p=1$ have nontrivial fixed point sets and other elements of
$\tilde T_p$ act freely. These fixed point sets are the same for all
$p$-roots of the unit and they are diffeomorphic to the $3$-dimensional
equator sphere in $S^5$. Let consider the ${\bold Z}_p$-action on $S^5$
which is given by the subgroup of $\tilde T_p$ formed by elements with
$z^p=1$. The factor-space of
$U(3)/SU(2)\cdot {\textstyle diag}(\bar w,\bar w,1)$ under this
${\bold Z}_p$-action is diffeomorphic to $S^5$ (one can consider
$S^5$ as the cyclic
$p$-covering of $S^5$ ramified at the $3$-dimensional
equator sphere). The factor-group $\tilde
T_p/{\bold Z}_p$ acts freely on this factor-space
and for every $p$ we obtain the mapping
$$
W^7_p \longrightarrow \left( \matrix z^{p} &
0 & 0 \\ 0 & z^{p} & 0 \\ 0 & 0 & z \endmatrix \right) {\big \backslash
U(3) /} SU(2) \cdot \left( \matrix {\bar w} & 0 & 0 \\ 0 & {\bar w} & 0
\\ 0 & 0 & 1 \endmatrix \right) = {\bold C}P^2.
$$

It is almost evident that this mapping is the ${\bold
R}P^3$-bundle but the strong proof needs an additional work.

Recently Ya. V. Bazaikin answering to our question computed that the
order of $H^4(W^7_p)$ is equal to $r_p=(4p-1)$.

The $13$-dimensional case is absolutely analogous to the $7$-dimensional
case and for every $p$ we obtain the mapping
$$
M^{13}_p \longrightarrow
T^1_p \backslash U(5) / SU(4) \cdot T_{(0)} = {\bold C}P^4.
$$

This mapping also ought to be an ${\bold R}P^5$-bundle.

It was computed in \cite{Ba} that the order of
$H^6(M^{13}_p)$ is equal to $s_p=(8p^2-4p+1)$.

One can notice that the very interesting formula
$$
s_p = \frac{r_p^2+1}{2}
$$
holds. For $p=1$ it coincides with (8) and (12).

One can also ask how to extend these
embeddings $W^7_p \rightarrow M^{13}_p$ to towers.

2) When we discussed the results of this paper with K. Grove he asked
about the existence of such embeddings for the even-dimensional Wallach
spaces (\cite{W}). The existence of topological embeddings of these flag
spaces is evident. We can give a simplest example of an involution
$$
A \rightarrow {\hat S} \cdot A \cdot {\hat S}^{-1},
$$
$$
A =
\left(
\matrix \sqrt{-1} I_3 & 0 \\ 0 & -\sqrt{-1} I_3
\endmatrix
\right).
$$
This involution generates an involution on the
space $Sp(3)/(Sp(1) \times Sp(1) \times Sp(1))$. The
component of the fixed point set, which contains the orbit of the unit,
is diffeomorphic to the space $SU(3)/T^2$ where $T^2$ is the maximal
torus. Thus we obtain the totally geodesic embedding $SU(3)/T^2
\rightarrow Sp(3)/(Sp(1) \times Sp(1) \times Sp(1))$.

We remind that in \cite{V} it was proved that homogeneous metrics
on the even-dimensional Wallach spaces have the same maximal
pinching which is equal to $1/64$.
Thus one can expect that  this embedding is
pinching-essential and that for the other pair of the Wallach spaces
$(F_4/Spin(8),Sp(3)/(Sp(1) \times Sp(1) \times Sp(1)))$ such
embedding also exists.

\enddocument